\RequirePackage{fix-cm}
\documentclass[smallextended]{svjour3}       % onecolumn (second format)
\smartqed  % flush right qed marks, e.g. at end of proof
\usepackage{graphicx}
\usepackage{amsfonts,bm}
\usepackage{amsmath}
\usepackage{color}
\begin{document}

\title{Criterion of positivity for semilinear problems with applications in biology
%\thanks{Grants or other notes
%about the article that should go on the front page should be
%placed here. General acknowledgments should be placed at the end of the article.}
}
%\subtitle{Do you have a subtitle?\\ If so, write it here}

%\titlerunning{Short form of title}        % if too long for running head

\author{Michel Duprez         \and
        Antoine Perasso %etc.
}

%\authorrunning{Short form of author list} % if too long for running head

\institute{M. Duprez \at
              Institut de Mathématiques de Marseille UMR CNRS 7373 \\
              Universit\'{e} Aix-Marseille\\
              39, rue F. Joliot Curie, 13453 Marseille Cedex 13, France\\
              Tel.: +33 (0) 4 13 55 14 65 \\
              %Fax: +123-45-678910\\
              \email{mduprez@math.cnrs.fr}\\
              Corresponding author           %  \\
%             \emph{Present address:} of F. Author  %  if needed
           \and
           A. Perasso \at
              Chrono-environnement UMR CNRS 6249\\
              Universit\'{e} Bourgogne Franche-Comt\'{e}\\
              16 route de Gray, 25000 Besan\c{c}on, France\\
              %Tel.: +33 (0)3 81 66 20 14 \\              
              %\email{antoine.perasso@univ-fcomte.fr}
}

\date{Received: date / Accepted: date}
% The correct dates will be entered by the editor

\maketitle

\begin{abstract}
The goal of this article is to provide an useful criterion of positivity and well-posedness for a wide range of infinite dimensional semilinear abstract Cauchy problems. This criterion is based on some weak assumptions on the non-linear part of the semilinear problem and on the existence of a strongly continuous semigroup generated by the differential operator. To illustrate a large variety of applications, we exhibit the feasibility of this criterion through three examples in mathematical biology: epidemiology, predator-prey interactions and oncology.
%Insert your abstract here. Include keywords, PACS and mathematical
%subject classification numbers as needed.
\keywords{Positivity \and Well-Posedness \and Dynamic Systems \and Semilinear Problems \and Population Dynamics}
% \PACS{PACS code1 \and PACS code2 \and more}
 \subclass{35A01 \and 35B09 \and 35Q92 \and 92D25}
\end{abstract}

\section{Introduction}
In a wide range of mathematical modelling of natural phenomena, the quantities that are described through the mathematical system have to satisfy some positivity properties to ensure physical reality. For instance, when considering the evolution of matter quantities, such as in biology 
(or also physics \cite{alaa2008mathematical}, chemistry \cite{turing1952chemical},...), the positivity of the solutions of the underlying dynamical system is a crucial prerequisite to achieve the well-posedness of the problem and to guarantee its physical relevance.

A significant proportion of dynamical systems that describe the evolution over time of matter quantities are non-linear, but it oftenly appears that the non-linear effects can be seen as perturbations of linear dynamics, leading to such a differential formulation: 
\begin{equation}\label{eq:genericmodel}
\begin{cases}
y'(t) = \overbrace{Ay(t)}^{\text{linear dynamics}} + \overbrace{f(y(t),t)}^{\text{perturbations}} , t\geq0,\\\noalign{\smallskip}
y(0) = y_0,
\end{cases}
\end{equation}
where $y(t)$ denotes the modeled matter quantity at time $t$, that mathematically lies in a Banach lattice. When imposing a non-negative initial condition $y_0$, the question of positivity is then crucial to study. In the case of a finite dimensional operator $A$, this question has been extensively studied (see  \cite{Smith95} and references therein for general results). However, to our knowledge, we don't know any general criterion of positivity in the case where $A$ is a differential operator, \textit{i.e.} when the first equality in \eqref{eq:genericmodel} rewrites as partial differential equations (PDEs), while such differential operators are extensively used in mathematical biology, or also in many other applied mathematical sciences. For instance, in the specific case of biology, let us mention the use of structured population dynamics models, where the operator is of transport type, or the use of diffusive processes, where models incorporate a Laplacian operator (see \cite{Pierre10} for a review of positivity results in reaction-diffusion systems).

The goal of this article is to provide an useful criterion of well-posedness and positivity for the semilinear problem \eqref{eq:genericmodel} for wide ranges of linear differential operators $A$ and non-linear functions $f$, and then to illustrate the feasibility of this criterion through three examples of models arising from mathematical biology: epidemiology, predation and oncology.

This article is structured as follows: Section \ref{sec:biol} is dedicated to the introduction of three concrete biological models, described by semilinear PDEs, for which the positivity of solutions must necessarilly be satisfied. Then we tackle in Section \ref{sec:criterion} the formulation and the proof of the criterion of positivity and well-posedness. This criterion is based on the formulation of an abstract semilinear Cauchy Problem, studied using a semigroup approach. Finally, in Section \ref{sec:Appl}, we apply the criterion to the biological models of Section \ref{sec:biol} to prove the well-posedness and the positivity of their solution.
%Section \ref{sec:Rem} is finally dedicated to final remarks and perspectives.
%\bigskip

%\textcolor{red}{Ce qui suit : à déplacer et/ou supprimer}
%\\
%
%The positivity of solutions in partial differential theory is crucial 
%in a wide variety of problems. 
% For example, when a dynamic system models the evolution of matter quantities, 
%it is very important that the solution of this model remains positive. 
%For a non-linear Cauchy problems, we will give a criterion on the differential operator 
%and the non-linearity for existence of positive solutions. 
%Moreover, we will work in a Banach lattice. Thus the results holds true 
%for  a lot in function spaces : $L^p(\Omega\times (0,T))$ 
%($\Omega\subset \mathbb{R}^N$, $N\in \mathbb{N}^*$, $p\in \mathbb{N}^*$), $\mathcal{C}(\overline{\Omega}\times [0,T])$...

\section{Three biological examples}\label{sec:biol}

In this section, we introduce three examples of semilinear evolutionary problems in mathematical biology for which the positivity and well-posedness have to be proved for biological purpose.  
The matter quantities that are modelled in those three examples, \textit{i.e.} populations, predator/prey or cell densities, evolve with respect to the time $t\geq 0$. The epidemiological and predator-prey models deal with transport process, with a non-constant velocity in the epidemiological case and a non-local boundary condition in the predator-prey case, while the model in oncology deals with diffusive PDEs.\\
One can note that through those specific examples, a large spectrum of biological models are involved: PDE structured population models (see \cite{Magal08} and references therein) and reaction-diffusion models.

\paragraph{Epidemiology}

The first example on which we focus deals with epidemiology. When modeling the transmission of disease between individuals, a common way is to split the population densities into two sub-classes that are the susceptible class ($S$) and the infected class ($I$). From such a splitting results the classical epidemiological model of SI type \cite{Kermack27}. Furthermore, lots of diseases (influenza, HIV, prion pathologies...) have a varying intensity during their evolution that may be important to take into account in the modeling process. This phenomena was recently described in \cite{perasso2013infection,perasso_raza14}, where the disease intensity was incorporated into the infected class, leading to the formulation of the following infection load-structured epidemiological model of transport type:
\begin{equation}\label{systcompletepidemio}
 \left\{\begin{array}{l}
S'(t)
=\gamma -(\mu_0+\alpha)S(t)- S(t)\mathcal{T}(\beta I)(t),~t\geq 0 ,\\\noalign{\smallskip}
\partial_t I(t,i)
=-\partial_i(\nu i I(t,i))-\mu (i)I(t,i)
	 +\phi(i)S(t)\mathcal{T}(\beta I)(t),~t\geq 0,~i\in J,\\
\noalign{\smallskip}\nu \kappa I(t,\kappa )=\alpha S(t),~t\geq 0 ,\\
\noalign{\smallskip} S(0) = S_0, \quad I(0,\cdot)=I_0(\cdot), \\
        \end{array}
\right.
\end{equation}
where the infection load is $i\in J=(\kappa ,+\infty)\subset \mathbb{R}_+$, $\mathcal{T}$ is the integral operator 
defined for some integrable fonction $h$ on $J$ by
\begin{equation*}
 \mathcal{T}:h\rightarrow \displaystyle\int_J h(i)di
\end{equation*}
and the epidemiological parameters satisfy the following assumptions:
\begin{itemize}%[(i)]
%  \item[-] $S_0\in\mathbb{R}_+$ and $I_0$ is a positive real valued function	 on $J$
  \item[-] $\beta, ~\mu_0,~\nu,~\alpha >0$ and $\gamma \geq 0$;
  \item[-] $\phi\in \mathcal{C}^{\infty}(J)$ is a non-negative function such that 
  $\lim\limits_{i\rightarrow +\infty}\phi(i)=0$ and $\int_J\phi(i)di=1$, $\mu\in L^{\infty}(J)$ is such that $\mu(i)\geq \mu_0$ for almost every $i\in J$.
 \end{itemize}
For a biological relevance, it is clear that for each positive initial condition $(S_0,I_0(\cdot))$, the densities $S(t)$ and $I(t,\cdot)$ in Problem (\ref{systcompletepidemio}) have to remain positive whenever they exist.

\paragraph{Predator-prey interactions} When considering predator-prey interactions, the age of the prey is a key factor of selection for the predator. It is therefore natural to add a structuration of the prey densities according to their age. In doing so, the classical Lotka-Volterra model, that was initially an ODE model \cite{Murray04}, turns into the following PDE model, that is developed in \cite{Perasso17}:
\begin{equation}\label{syst preda}
\left\{\begin{array}{l}
\partial_t x(t,a)+\partial_a x(t,a)=-\mu(a) x(t,a)- y(t) \gamma(a) x(t,a),
~t\geq 0,~a\geq0,\\\noalign{\smallskip}
 y'(t)=\alpha y(t) \int_0^{\infty}\gamma(a)x(t,a)da-\delta y(t),~t\geq 0,~a\geq0,\\\noalign{\smallskip}
 x(t,0)=\int_0^{\infty}\beta(a)x(t,a)da,~t\geq0\\\noalign{\smallskip}
x(0,\cdot)=x_0(\cdot),~y(0)=y_0,
\end{array}\right.
\end{equation}
where $x$ and $y$ denote the density of preys and predators, respectively. The assumptions on the parameters are the following:
\begin{itemize}
\item[-] $\alpha \in ]0,1[$, $\delta>0$ are constant parameters that respectively denote the assimilation coefficient of ingested preys and the basic mortality rate of the predators;
\item[-] $\mu,\gamma,\beta\in L_+^\infty(\mathbb{R}_+)$ are age-dependent functions that represent, 
respectively, the basic mortality rate of the preys, the predation rate and the birth rate.
\end{itemize}
To ensure a certain realism, we want that the densities of preys $x$ and predators $y$ remain positive given a positive initial data $(x_0,y_0)$.

\paragraph{Oncology} The third application is a model that describes the growth of a brain tumour published in \cite{chakrabarty2009distributed}. The model aims at studying a treatment method of tumor cells through a problem of controllability. The tumor and normal cells are in competition  for the resources and are subject to a drug treatment whose role is to decrease the cell densities. Even if some normal cells are destroyed, the key point here is that the drug affects more the tumor ones.\\
To make explicit the model, let us consider $\Omega$ a bounded domain of $\mathbb{R}^N$, $N\in \mathbb{N}^*$, 
with boundary $\partial \Omega$ of class $\mathcal{C}^2$ and for a fixed $T>0$, let $Q_T=\Omega\times (0,T)$ 
and $\Sigma_T=\partial \Omega\times (0,T)$. 
The evolution problem is then written using the following three semilinear heat equations, where the variables $(t,x)$ are delibarately avoided for a better reading:
\begin{equation}\label{systcomplet}
 \left\{\begin{array}{l}
\partial_t y_1=d_1\Delta y_1
+a_1g_1(y_1)y_1-(\alpha_{1,2}y_2+\kappa_{1,3}y_3)y_1,\\ \noalign{\smallskip}
\partial_t y_2=d_2\Delta y_2
+a_2g_2(y_2)y_2-(\alpha_{2,1}y_1+\kappa_{2,3}y_3)y_2,\\\noalign{\smallskip}
\partial_t y_3=d_3\Delta y_3 - a_3 y_3+u ,\\\noalign{\smallskip}
\partial_ny_i(t)=\nabla y_i(t)\cdot\overrightarrow{n}
=0,~t\geq 0, i\in\{1,...,3\},\\\noalign{\smallskip}
y(x,0)=y_0(x),~x\in\Omega,
        \end{array}
\right.
\end{equation}
where  $y:=(y_1,y_2,y_3)^*$, 
$\overrightarrow{n}$  denotes the external normalized normal to the boundary $\partial \Omega$. Here $y_1(t,x)$ stands for the density of tumor cells, $y_2(t,x)$ the density of normal tissue and $y_3(t,x)$ the drug concentration at any vector position $x$ and time $t$. In the latter problem, the growth rates of cells are defined by the functions $g_i$ according to the following logistic shape:
\begin{equation*}
 g_i(y_i)=1-y_i/k_i.
\end{equation*}
The assumptions on the parameters are the following:
\begin{itemize}
\item[-] $d_i>0$ are the coefficients for the space diffusive effect;
\item[-]  $a_i>0$,  where $a_1$, resp. $a_2$, denotes the tumor
cell intrinsic growth rate, resp. the normal tissue intrinsic growth rate and $a_3$ is the drug reabsorption coefficient;
\item[-]  $k_i>0$ denote the carraying capacity of the medium;
\item[-]  $\alpha_{i,j}>0$ are coefficients that translate the interspecific competition between tumor and normal cells;
\item[-]  $\kappa_{1,3}\gg\kappa_{2,3}>0$ are the degradation rates due to the treatment;
\item[-]  $u(x,t)\geq 0$ represents the flux of injected drug over time at position $x$.
\end{itemize}
%Note that the term $-a_3 y_3$, called reabsorption effect of the drug, is due to the chemical reaction of drug with tissues.\\
Similarly to the previous biological examples, we aim at proving well-posedness and positivity of the solution.

\section{A criterion of positivity and well-posedness}\label{sec:criterion}

%\subsection{Setting}

In all this section, let us consider $(\mathcal{W},+,\|\cdot\|_{\mathcal{W}},\geq)$ a Banach lattice (see \cite[p. 6]{meyer1991banach}), 
\textit{i.e.} an partially ordered Banach space for which any given elements $x,y$ of  $\mathcal{W}$ 
have a supremum $\sup(x,y)$ %and an infimum $\inf(x,y)$
and %where the norme and the order are compatible, i.e.  
for all $y_1,~y_2,~y_3\in \mathcal{W}$ and $\alpha\geq0$,
\begin{equation}\label{prop banach lattice}
\left\{\begin{array}{l}
y_1\leq y_2\Rightarrow (y_1+y_3\leq y_2+y_3\mbox{ and }\alpha y_1\leq\alpha y_2),\\\noalign{\smallskip}
 |y_1|_{\mathcal{W}}\leq  |y_2|_{\mathcal{W}} 
 \Rightarrow  \|y_1\|_{\mathcal{W}}\leq  \|y_2\|_{\mathcal{W}},
\end{array}\right.
\end{equation}
with, for all $y\in\mathcal{W}$, $|y|_{\mathcal{W}}=\sup(y,\,-y)$.  
We will denote by $\mathcal{W}^+=\{y\in \mathcal{W} : 0\leq y\}$ the non-negative cone and for every $m>0$ by $B_m$ the ball of $\mathcal W$ of radius $m$.

We consider in this work the system
\begin{equation}\label{equation nonlineraire generale pa2}
\left\{\begin{array}{ll}
y'(t)=Ay(t)+f(y(t),t),~t\geq0&\mathrm{in}~\mathcal{W},\\\noalign{\smallskip}
y(0)=y_0&\mathrm{in}~\mathcal{W},%\mathrm{in}~\Omega,
\end{array}\right.
\end{equation}
where $A:D(A)\subset\mathcal{W}\rightarrow\mathcal{W}$ is an infinitesimal generator 
of a  $\mathcal{C}_0$-semigroup 
$(T_A(t))_{t\geq0}$, $y'(t)$ is an element of $\mathcal{W}$ and  $f:\mathcal{W}\times\mathbb{R}^+\rightarrow\mathcal{W}$  is continuous in $t$ and locally Lipschitz continuous in $y$ uniformly in $t$ in the following sense: for every $m>0$ there exists a constant $k_m>0$ such that for every $z_1,z_2\in B_m$,
\begin{equation*}
\|f(z_1,t)-f(z_2,t)\|_{\mathcal{W}}\leq  k_m \|z_1-z_2\|_{\mathcal{W}},~\forall  t\in\mathbb{R}^+.
\end{equation*}

Finally, let us briefly remind that for a fixed $T\in]0,\infty]$, a \textit{mild solution} of Problem \eqref{equation nonlineraire generale pa2} on $[0,T[$ is a function $y\in\mathcal{C}([0,T[;\mathcal{W})$ that satifies the integral equation
\begin{equation*}
 y(t)=T_A(t)y_0+\displaystyle\int^t_{0}T_A(t-s)f(y(s),s)\mathrm{d}s.
\end{equation*}
%And a function $y:[0,T[\rightarrow\mathcal{W}$ will be a \textit{classical solution} 
%of the initial value problem (\ref{equation nonlineraire generale pa2}) 
%on $[0,T[$ if $y$ is continuous on $[0,T[$, continuously differentiable 
%on $]0,T[$, $y(t)\in D(A)$ for all $t\in ]0,T[$ and (\ref{equation nonlineraire generale pa2}) 
%is satisfied on $[0,T[$.

\begin{remark}
 Since $\mathcal{W}^+$ is closed (see \cite{meyer1991banach}), 
 we deduce that for all $T>0$, the order $\geq$ is compatible 
 with the integration in time, more precisely,  
for all $x, y\in \mathcal{C}([0,T];\mathcal{W})$,
\begin{equation}\label{compatibilite integral}
\left(x(t)\geq y(t)~ \forall~ t\in[0,T]\right)
\Rightarrow \displaystyle\int^T_0 x(s)\mathrm{d}s\geq 
\displaystyle\int^T_0 y(s)\mathrm{d}s.
\end{equation}
\end{remark}

The following theorem, that states well-posedness and positivity property for the solution of Problem \eqref{equation nonlineraire generale pa2}, is the main result of the present article:
\begin{theorem}\label{theo exist local sol pa3}
Let  $y_0\in \mathcal{W}^+$. We suppose that 
\begin{enumerate}
\item[(i)] $A$ is generator of a positive $\mathcal{C}_0$-semigroup on $\mathcal{W}$,
\textit{i.e.} $T_A(t)\mathcal{W}^+\subset \mathcal{W}^+$ for all $t\geq 0$,
%\item[(ii)] $f$ is locally Lipschitz continuous in $y$ and continuous in $t$,
\item[(ii)] for all $m>0$, there exists $\lambda_m\in\mathbb{R}$ 
such that, for all $z\in \mathcal{C}(\mathbb{R}^+; \mathcal{W}^+\cap B(0,m))$, 
\begin{equation}\label{cond positive A+B2*}
 f(z(t),t)+\lambda_m z(t)
\geq 0, \quad \forall t\geq 0.
\end{equation}
\end{enumerate}
Then there exists $t_{max}\in]0,\infty]$ such that system \eqref{equation nonlineraire generale pa2} 
has an unique positive mild solution $y\in\mathcal{C}([0,t_{max}[; \mathcal{W})$. 
Moreover, if $t_{max}<\infty$, 
\begin{equation*}
\lim\limits_{t\rightarrow t_{max}}\|y(t)\|_{\mathcal{W}}=\infty.
\end{equation*} 
\end{theorem}

The main idea of the proof is to perform a vectorial translation to the range values of the non-linear part $f$ so that they remain in $\mathcal W^+$. This translation is then compensated by the substraction of a linear term to the differential operator, that does not affect its spectral and positivity properties. Consequently, we shall study the following system in the proof of the theorem:
\begin{equation}\label{equation modifie B}
\left\{\begin{array}{ll}
 y'(t)=(A-\lambda I)y(t)
 +f(y(t),t)+\lambda y(t),~t>0&\mathrm{in}~\mathcal{W},\\\noalign{\smallskip}
y(0)=y_0&\mathrm{in}~\mathcal{W}.
\end{array}\right.
\end{equation}

\begin{remark}\label{rmk 2}
Since $A$ is an infinitesimal generator of a positive $\mathcal{C}_0$-semigroup $(T_A(t))_{t\geq0}$, 
then, for every $\lambda\in\mathbb{R}$, $A-\lambda I$ is also an infinitesimal generator of a positive $\mathcal{C}_0$-semigroup  $(T_{A-\lambda I}(t))_{t\geq0}$.
Indeed, we remark that $T_{A-\lambda I}(t)=e^{-\lambda t}T_A(t)$ for all $t\geq 0$.
\end{remark}

%\begin{lemma}\label{lemme}
%Let $A$ be an infinitesimal generator of a positive $\mathcal{C}_0$-semigroup. Then, for every $\lambda\in\mathbb{R}$, $A-\lambda I$ is an infinitesimal generator of a positive $\mathcal{C}_0$-semigroup.
%\end{lemma}
%
%\begin{proof}
%As a bounded perturbation of $A$, $A-\lambda I$ is an infinitesimal generator 
%of a $\mathcal{C}_0$-semigroup $(T_{A-\lambda I}(t))_{t\geq 0}$ on $\mathcal{W}$ 
%(see \cite[p. 76]{pazy1983semigroups}). 
%A $\mathcal{C}_0$-semigroup on a Banach lattice is positive if and only if 
%the resolvent $(\mu I-L)^{-1}$ of its generator $L$ is positive 
%for all sufficiently large $\mu$ 
%(see \cite[p. 207]{engel2006short}). 
%Thus there exists $\mu^*$ such that, for all $x\in \mathcal{W}$ and all $\mu>\mu^*$,
%$ (\mu I-A)^{-1}x\geq 0$.
%Consequently, for all $x\in \mathcal{W}$ and  all $\mu>\mu^*-\lambda$, we have
%$ (\mu I-A+\lambda I)^{-1}x\geq 0$.
%Then $A-\lambda I$ is an infinitesimal generator of a positive 
%$\mathcal{C}_0$-semigroup. 
%\end{proof}

\begin{proof}[Proof of Theorem 1.1]
Without loss of generality, 
we can assume that $\lambda_m$  is nonnegative in \eqref{cond positive A+B2*}.
Since $A$ is generator of a positive $\mathcal{C}_0$-semigroup $(T_{A}(t))_{t\geq 0}$, there exists $\omega, M\geq 1$ such that, for all $t\in\mathbb{R}^+$,
\begin{equation*}
 \|T_{A}(t)\|_{\mathcal{W}}\leq Me^{\omega t}.
\end{equation*}
Remark \ref{rmk 2} then implies that for evey $\lambda\in\mathbb{R}$, $A-\lambda I$ is also generator of a positive $\mathcal{C}_0$-semigroup $(T_{A-\lambda I}(t))_{t\geq 0}$. Moreover, it is easy to check that for all $t\in\mathbb{R}^+$,
\begin{equation}\label{estim semigroup}
 \|T_{A-\lambda I}(t)\|_{\mathcal{W}}\leq Me^{\omega t}, \quad \forall \lambda\in\mathbb{R}^+.
\end{equation}

Let $t_0\in(0,1)$, $m=2 M \, e^{\omega}\|y_0\|_{\mathcal{W}}$ and $\lambda_m$ that satisfies \eqref{cond positive A+B2*}. Consider the set $\Gamma_m =\{y\in \mathcal{C}([0,t_0];\mathcal{W}^+):y(0)=y_0, \|y(t)\|_{\mathcal{W}}\leq m,\forall t\in[0,t_0]\}$. The continuity properties of the lattice operations (see \cite{meyer1991banach}) imply that $\Gamma_m$ is a non-empty closed subset of $\mathcal{C}([0,t_0];\mathcal{W})$.

Consider now the mapping $\psi$, defined on $\Gamma_m$ by
\begin{equation*}
  \psi(y)(t)=T_{A-\lambda_m I}(t)y_0+\displaystyle\int_{0}^{t}T_{A-\lambda I}(t-s)
 \left[f(y(s),s)+\lambda_m y(s)\right]\mathrm{d}s, \quad  t\in[0,t_0].
\end{equation*}
We aim at proving that $\psi$ has a unique fixed point in $\Gamma_m$.\\
Let us start by proving that $\psi$ preserves $\Gamma_m$. The positivity of $(T_{A-\lambda_m I}(t))_{t\geq 0}$ and the positivity assumption \eqref{cond positive A+B2*} clearly imply that $\psi(y) \in \mathcal{C}([0,t_0];\mathcal{W}^+)$. Furthermore, from the inequality \eqref{estim semigroup},  one deduces that
\begin{equation*}
\begin{array}{l}
\|\psi(y)(t)\|_{\mathcal{W}} \leq 
Me^{\omega t}\|y_0\|_{\mathcal{W}}
+Me^{\omega t}\displaystyle\int_{0}^{t_0}(\|f(y(s),s)-f(0,s)\|_{\mathcal{W}}\\
\hspace*{6cm}+\|f(0,s)\|_{\mathcal{W}}+\lambda_m\|y(s)\|_{\mathcal{W}})\mathrm{d}s.
\end{array}
\end{equation*}
The time continuity property on $f$ induces the existence of $\gamma>0$ (independent of $t_0<1$) such that for every $y\in\Gamma_m$ and every $t\in(0,t_0)$,
\begin{equation*}
\|\psi(y)(t)\|_{\mathcal{W}} \leq
Me^{\omega}(\|y_0\|_{\mathcal{W}}+t_0(mk_m +\gamma+m\lambda_m )).
\end{equation*}
Thus, for $t_0=\min\{1,\|y_0\|_{\mathcal{W}}\times(mk_m +\gamma+m\lambda_m )^{-1}\}$ we have $\|\psi(y)(t)\|_\mathcal{W} \leq 2 M e^{\omega} \|y_0\|_\mathcal{W} = m$ and so $\psi(y)\in\Gamma_m$. 
\\
We now prove that $\psi$ is contractant in the following sense: for every $y,z\in \Gamma_m$, every $n\in\mathbb{N}^*$ and every $t\in[0,t_0]$,
\begin{equation}\label{recurence exist glob bornee}
\|\psi^n(y)(t)-\psi^n(z)(t)\|_{\mathcal{W}}
\leq  \dfrac{[Me^{\omega }t(k_m+\lambda_m)]^n}{n!}
\sup\limits_{t\in \left[0,t_0\right]}\|y(t)-z(t)\|_{\mathcal{W}}.
\end{equation}
Let us prove \eqref{recurence exist glob bornee} by induction. 
By definition of $\Gamma_m$, we have $$\|y(t)\|_{\mathcal{W}},~\|z(t)\|_{\mathcal{W}}\leq m$$
for all $t\in[0,t_0]$. 
Then the Lipschitz assumption on $f$ implies that
\begin{equation*}
\|\psi(y)(t)-\psi(z)(t)\|_{\mathcal{W}} \leq Me^{\omega }\left(k_m+\lambda_m\right)t
\sup\limits_{\theta\in \left[0,t_0\right]}\|y(\theta)
-z(\theta)\|_{\mathcal{W}},
\end{equation*}
and equality \eqref{recurence exist glob bornee} holds for $n=1$. 
Suppose now that \eqref{recurence exist glob bornee} holds for a $k\in\mathbb{N}^*$. Then for all $t\in [0,t_0]$,
\begin{equation*}\begin{array}{rl}
%\sup\limits_{t\in \left[0,T\right]}
&\|\psi^{k+1}(y)(t)-\psi^{k+1}(z)(t)\|_{\mathcal{W}}\\
&\leq  (Me^{\omega }(k_m+\lambda_m))
\displaystyle\int_{0}^{t}\|\psi^k(y)(s)-\psi^k(z)(s)\|_{\mathcal{W}}\mathrm{d}s,\\\noalign{\smallskip}
&\leq \dfrac{[Me^{\omega }(k_m+\lambda_m)]^{k+1}}{k!}
\sup\limits_{\theta\in \left[0,t_0\right]}\|y(\theta)-z(\theta)\|_{\mathcal{W}}\displaystyle\int_{0}^{t}s^k\mathrm{d}s,
\end{array}
\end{equation*}
and \eqref{recurence exist glob bornee} is true for $k+1$ and consequently for every $n\in\mathbb{N}^*$ by induction. Finally, we can apply the Banach's fixed point theorem 
to conclude that $\psi$  has a unique fixed point $\bar{y}$ in $\Gamma_m$.  Systems \eqref{equation nonlineraire generale pa2} and \eqref{equation modifie B} being equivalent, $\bar y$ is a mild solution of \eqref{equation nonlineraire generale pa2}. Then some standard time extending properties of the solution induce that the solution $\bar y$ is defined on a maximal interval $[0,t_{max}[$. 
To finish, we prove the uniqueness of the solution on the whole space 
$\mathcal{C}([0,t_{max}(\bar y)[,\mathcal{W}^+)$.
If $\bar z$ is another mild solution defined on $[0,t_1[$ with $t_1<t_{max}(\bar y)$, then, denoting $R=\max\limits_{\theta\in [0,t_1]}\{\|\bar y(\theta)\|_{\mathcal{W}},\|\bar z(\theta)\|_{\mathcal{W}}\}$, we obtain for all $t\in[0,t_1]$,
\begin{equation*}
 \|\bar y(t)- \bar z(t)\|_{\mathcal{W}} \leq Me^{\omega t_1}k_R
\displaystyle\int_{0}^{t}\|\bar y(s)- \bar z(s)\|_{\mathcal{W}}\mathrm{d}s.
\end{equation*}
Then $\|\bar y(t)- \bar z(t)\|_{\mathcal{W}}= 0$ by a standard Gronwall argument and 
$\bar y=\bar z$ in $[0,t_1]\times \mathcal{W}$. 
Furthermore, if $t_{max}(\bar y)<\infty$, since $\|\bar z(t)\|_{\mathcal{W}}=\|\bar y(t)\|_{\mathcal{W}}$ for all $t<\min\{t_{max}(\bar y),t_{max}(\bar z)\}$ 
and $\lim\limits_{t\rightarrow t_{max}(\bar y)}\|\bar y(t)\|_{\mathcal{W}}=\infty$, 
we deduce that the maximal intervals of existence of $\bar y$ and $\bar z$ are equal.

% We have just prove  there exists a maximal  $T^*$ such that we have a unique mild solution in $\mathcal{C}([0,T^*[;\mathcal{W})$.% and it is positive. 
% Moreover if $T^*<\infty$, then $\lim\limits_{t\rightarrow T^*}\|y(t)\|_{\mathcal{W}}=\infty$. Indeed otherwise there exists a sequence 
% $t_n\rightarrow T^*$ such that
% \begin{equation*}
%  \|y(t_n)\|_{\mathcal{V}}\leq C.
% \end{equation*}
% In the preview argument we can chose $M_{\lambda_m}:=Re^{\omega(T^*-\lambda_m)}$. Then 
% %If we set $r=\min\{n:T^*-t_n\leq1\}$, we can chose in the precedent argument for all $n\geq r$ $M(t_n)=M(t_r)$ and 
% %$\tilde{y}_{max}(t_n)=\tilde{y}_{max}(t_r)$. In this way, 
% $\delta$ is independant of $t_n$ and we obtain a contradiction with the definition of $T^*$.
\end{proof}

\section{Illustrations of the criterion in mathematical biology}\label{sec:Appl}

In this section, we exhibit the application of well-posedness and positivity criterion on the three biological examples of Section \ref{sec:biol}.

\paragraph{Epidemiology}

Consider the Banach lattice  $X=\mathbb{R}\times L^1(J)$,
$X^+$  the non-negative cone of $X$ and $y_0=(S_0,I_0)\in X^+$.
%Using Theorem \ref{theo exist local sol pa3} we can obtain the following result:
%
%\begin{theo}
% For every $(S_0,I_0)\in X^+$,  problem (\ref{systcompletepidemio}) 
% has an unique mild solution $(S,I)\in \mathcal{C}(\mathbb{R}^+,X^+)$. 
%\end{theo}
Then it is clear that Problem \eqref{systcompletepidemio} can rewrite as \eqref{equation nonlineraire generale pa2}, 
where the function $f:X\rightarrow X$ and 
the differential operator $A:D(A)\subset X\rightarrow X$ 
are given by
\begin{equation*}
\begin{array}{c}
 f(u,v)=
 \left(\begin{array}{c}f_1(u,v)\\
f_2(u,v) \end{array}\right)
=\left(\begin{array}{c}\gamma-u\mathcal{T}(\beta v)\\
\phi u\mathcal{T}(\beta v) \end{array}\right),~
 A=\left(\begin{array}{cc}
        -\mu_0  -\alpha&0\\
        0&-\frac{d}{di}(\nu i\cdot )-\mu
         \end{array}
\right),\\\noalign{\smallskip}

\end{array}\end{equation*}
with $D(A)=\{(x,\varphi)\in X, (i\varphi)\in W^1_1(J)\mathrm{~and~}\nu \kappa\varphi (\kappa ) =\alpha x\}$. 
%$ L_1=-\mu_0  -\alpha$ and  $L_2=-\frac{d}{di}(\nu i\cdot )-\mu$.
%\begin{equation*}
%\left\{\begin{array}{l} f_1(u,v)=\gamma-u\mathcal{T}(\beta v)\\
% f_2(u,v)=\phi u\mathcal{T}(\beta v)
%\end{array}\right.
%\mbox{ and }
%\left\{\begin{array}{l} L_1=-\mu_0  -\alpha.\\
% L_2\varphi=-\dfrac{d}{di}(\nu i\varphi)-\mu\varphi
%\end{array}\right.
%\end{equation*}
In \cite{perasso2013infection}, the authors prove that 
 the differential operator $(A,D(A))$ is an infinitesimal generator of 
 a  positive $\mathcal{C}_0$-semigroup $(T_A(t))_{t\geq 0}$ on X 
 and that function $f$ is locally Lipschitz continuous on $X$. 
Moreover, for every $m>0$ and every 
$(\bar S, \bar I)\in \mathcal{C}(\mathbb{R}^+; X^+\cap B(0,m))$, one gets, denoting $\lambda_m = m\beta$,
 \begin{equation*}
 \left\{\begin{array}{l}
  f_1(\bar S(t),\bar I(t))+\lambda_m \bar S(t)
\geq \gamma+ \bar S(t)(\lambda_m-\beta \mathcal{T}(\bar I(t))) 
\geq 0,\\
 f_2(\bar S(t),\bar I(t))+\lambda_m  \bar I(t)
=\phi \bar S(t)\mathcal{T}(\beta \bar I(t)) +\lambda_m  \bar I(t)
 \geq 0.
 \end{array}\right.
 \end{equation*}
Thus, condition \eqref{cond positive A+B2*} of Theorem \ref{theo exist local sol pa3} is satisfied and
there exists $t_{max}\in]0,\infty]$ such that Problem \eqref{systcompletepidemio}  
 has an unique mild solution $(S,I)$ in $\mathcal{C}([0,t_{max}[,X^+)$.

\paragraph{Predator-prey interactions}

Let $X=L^1(\mathbb{R}^+)\times \mathbb{R}$, $X^+$ the non-negative cone 
and  $(x_0,y_0)\in X^+$. 
Considering the operator $A:D(A)\subset X\rightarrow X$ and the functional  $f:X \to X$ given by
\begin{equation*}
\begin{array}{c}f(\phi,z)=
\left(\begin{array}{c}
f_1(\phi,z)\\
f_2(\phi,z)
\end{array}\right)=
\left(\begin{array}{c}
- z\gamma\phi\\
\alpha z\int^{\infty}_0\gamma(a)\phi(a)da
\end{array}\right), ~
 A=\left(\begin{array}{cc}
        L&0\\
        0&-\delta
         \end{array}
\right),
\end{array}
\end{equation*}
with $ D(A)=\{(\phi,z)\in X, \phi\in W^1_1(\mathbb{R}^+)\mathrm{~and~}\varphi (0)
 =\int_0^{\infty}\beta(a)\phi(a)da\}$ and $L\phi=-\phi'-\mu \phi$. 
The map $f$ is clearly locally Lipschitz continuous on $X$. Furthermore, under the assumption that there exists $\mu_0>0$ such that $\mu(a)\geq \mu_0$ f.a.e. $a\in\mathbb{R}^+$  the operator $A$ is the infinitesimal generator of a  positive $\mathcal{C}_0$-semigroup $(T_A(t))_{t\geq0}$ on $X$.
This is a standard result 
that we can find for example in \cite{nagel86}.
Then,  for all $m>0$, denoting $\lambda_m= m\gamma$,  we obtain for every
 $(\bar x, \bar y)\in \mathcal{C}(\mathbb{R}^+; X^+\cap B(0,m))$ 
 \begin{equation*}
 \left\{\begin{array}{l}
  f_1(\bar x(t), \bar y(t))+\lambda_m \bar x(t)\geq \bar x(t)(\lambda_m -\alpha m \gamma)\geq 0,\\
 f_2( \bar x(t),\bar y(t))+\lambda_m  \bar y(t)\geq  0.
 \end{array}\right.
 \end{equation*}
Again, condition \eqref{cond positive A+B2*} of Theorem \ref{theo exist local sol pa3} holds and the existence of  $t_{max}\in]0,\infty]$ such that  system \eqref{syst preda}  
 has an unique mild solution $(x,y)$ in $\mathcal{C}([0,t_{max}[,X^+)$ is ensured.

\paragraph{Oncology}

Let $X=L^2(\Omega;\mathbb{R}^3)$, $X^+$ the corresponding non-negative cone, 
$y_0\in X^+$ and $u\in L^2(Q_T)^+$. 
Then system \eqref{systcomplet} can be reformulated as \eqref{equation nonlineraire generale pa2} 
where
\begin{equation*}\left\{\begin{array}{l}
% D=\mathrm{diag}(d_1,d_2,d_3),\\\noalign{\smallskip}
 f(y)=(g+h)(y) +(0,0,u)^*,\\\noalign{\smallskip}
g(y)=\mathrm{diag}(a_1g_1(y_1)y_1,a_2g_2(y_2)y_2,-a_3y_3),\\\noalign{\smallskip}
h(y)=\mathrm{diag}(-(\alpha_{1,2}y_2+\kappa_{1,3}y_3)y_1,-(\alpha_{2,1}y_1+\kappa_{2,3}y_3)y_2,0),\\\noalign{\smallskip}
%\mbf{U}=(0,0,u)^*,\\%(u_1,u_2,u_3).
A=\mathrm{diag}(d_1\Delta,d_2\Delta,d_3\Delta).
\end{array}\right.
\end{equation*}
% and the the operator $ \mbf{A}$ defined by %(\ref{definition A neumann2}).
% \begin{equation*}%\label{definition A neumann2}
% \begin{array}{rcl}
%  \mbf{A}:%&\mathbb{H}^1(\Omega)&\rightarrow\mathbb{H}^1(\Omega)'\\
% &u&\mapsto \left(\varphi\mapsto \langle Au,\varphi\rangle_{\mathbb{H}^1(\Omega)',\mathbb{H}^1(\Omega)}
% =\langle \nabla u,\nabla\varphi\rangle_{\mathbb{L}^2(\Omega)}\right).
% \end{array}
% \end{equation*}
%\begin{Lemme}
% The operator $A:D(A)=\{z\in H^2(\Omega;\mathbb{R}^3)
% :\partial_{n}z=0\}
% \subset L^2(\Omega;\mathbb{R}^3)\rightarrow L^2(\Omega;\mathbb{R}^3)$ 
% is the infinitesimal generator of a positive $\mathcal{C}_0$ semigroup $(T_A(t))_{t\geq0}$ 
% on $L^2(\Omega;\mathbb{R}^3)$.
%%  $W(0,T)^3$ where 
%%  $W(0,T)=\{y\in L^2(0,T;H^1(\Omega);\dfrac{\partial y}{\partial t}
%%   \in  L^2(0,T;H^1(\Omega)')\}$.
%\end{Lemme}
%
%\begin{proof}
The existence of the semigroup $(T_A(t))_{t\geq0}$ is a consequence of the  Lumer-Phillips Theorem (see \cite[p. 14]{pazy1983semigroups}) for maximal dissipative operators. Indeed, in the present case,  $A$ is clearly maximal dissipative since it is defined with Laplacian operators.
Using the maximum principle of the heat equation, the semigroup is positive. 
Consequently, when taking 
$\lambda_m=\max\{m(a_1/k_1-\alpha_{1,2}-\kappa_{1,3}),m(a_2/k_2-\alpha_{2,1}-\kappa_{2,3}),a_3\}$ for $m>0$, we obtain the following estimations for all  
$\bar y=(\bar y_1, \bar y_2, \bar y_3)\in \mathcal{C}(\mathbb{R}^+; X^+\cap B(0,m))$ 
\begin{equation*}\left\{\begin{array}{rcl}
f_1(\bar y)+\lambda_m \bar y_1&=&
 a_1g_1(\bar y_1)\bar y_1-(\alpha_{1,2}\bar y_2+\kappa_{1,3}\bar y_3)\bar y_1
+\lambda_m \bar y_1\\
&\geq& \bar y_1[\lambda_m-m(a_1/k_1-\alpha_{1,2}-\kappa_{1,3})]\geq0,\\\noalign{\smallskip}
f_2(\bar y)+\lambda_m \bar y_2&=&
 a_2g_2(\bar y_2)\bar y_2-(\alpha_{2,1}\bar y_1+\kappa_{2,3}\bar y_3)\bar y_2
 +\lambda_m \bar y_2\\
 &\geq& \bar y_2[\lambda_m-m(a_2/k_2-\alpha_{2,1}-\kappa_{2,3})]\geq0,\\\noalign{\smallskip}
f_3(\bar y)+\lambda_m \bar y_3&=&
 -a_3\bar y_3+u+\lambda_m \bar y_3\geq \bar y_3(\lambda_m-a_3)\geq0.
\end{array}\right.
\end{equation*}
Thus condition \eqref{cond positive A+B2*} is satisfied 
and, using Theorem \ref{theo exist local sol pa3}, 
there exists $t_{max}\in]0,\infty]$ such that  problem \eqref{systcomplet} 
 has an unique mild solution $(x,y)$ in $\mathcal{C}([0,t_{max}[,X^+)$.


\begin{thebibliography}{00}
%%
%%%% \bibitem{label}
%%%% Text of bibliographic item
%%
%%\bibitem{}

\bibitem{alaa2008mathematical}
N.~Alaa, I.~Fatmi, J.-R. Roche, A.~Tounsi, Mathematical analysis for a model of
  nickel-iron alloy electrodeposition on rotating disk electrode: parabolic
  case, International Journal of Mathematics and Statistics 2 (2008) 30--49.
  
  
    \bibitem{nagel86}
W.~Arendt, A.~Grabosch, G.~Greiner, U.~Groh, H.~P.~Lotz, U.~Moustakas, R.~Nagel, F.~Neubrander, U.~Schlotterbeck, One-parameter semigroups of positive operators, Lect. Notes in Math.,
vol. 1184. Springer-Verlag, 1986.
  
  

  
  
  
  
  
%  \bibitem{engel2006short}
%W.~Arendt, A.~Grabosch, G.~Greiner, U.~Groh, H.~P.~Lotz, U.~Moustakas, R.~Nagel, F.~Neubrander, U.~Schlotterbeck, One-parameter semigroups of positive operators, Springer Science+
%  Business Media, 2006.
%  
%      
% \bibitem{nagel86}
%K.-J. Engel, R.~Nagel, A short course on operator semigroups, Lect. Notes in Math., 
%vol. 1184. Springer-Verlag, 1986.
  
  
  
  
  \bibitem{chakrabarty2009distributed}
S.~Chakrabarty, F.~B. Hanson, Distributed parameters deterministic model for
  treatment of brain tumors using galerkin finite element method, Math. biosci.
  219~(2) (2009) 129--141.

   \bibitem{engel2006short}
K.-J. Engel, R.~Nagel, A short course on operator semigroups, Springer Science+
  Business Media, 2006.

\bibitem{Kermack27}
W.~O. Kermack, M.~A. G., A contribution to the mathematical theory of
  epidemics, Proc. R. Soc. Lond. Ser. A 219 (1927) 700--721.

\bibitem{Magal08}
P.~Magal, S.~Ruan, Structured Population Models in Biology and Epidemiology,
  Vol. 1936 of Lecture Notes in Mathematics / Mathematical Biosciences
  Subseries, Springer, 2008.


\bibitem{meyer1991banach}
P.~Meyer-Nieberg, Banach lattices, Universitext, Springer-Verlag, Berlin, 1991.

\bibitem{Murray04}
J.~Murray, Mathematical Biology I, An introduction, third edition Edition,
  Interdisciplinary applied mathematics, Springer, 2004.
  
  

\bibitem{pazy1983semigroups}
A.~Pazy, Semigroups of linear operators and applications to partial
  differential equations, Vol.~44 of Applied Mathematical Sciences,
  Springer-Verlag, New York, 1983.


\bibitem{perasso2013infection}
A.~Perasso, U.~Razafison, Infection load structured si model with exponential
  velocity and external source of contamination, in: World Congress on
  Engineering, 2013, pp. 263--267.
  
  
\bibitem{perasso_raza14}
A.~Perasso, U.~Razafison, Asymptotic behavior and numerical simulations for an
  infection load-structured epidemiological model: application to the
  transmission of prion pathologies, SIAM J. Appl. Math. 74~(5) (2014)
  1571--1597.
  
  
  \bibitem{Perasso17}
A.~Perasso, Q.~Richard, Implication of age-structure on the dynamics of Lotka Volterra
equations, to appear in Differential and Integral Equations.

\bibitem{Pierre10}
M.~Pierre, Global existence in reaction-diffusion systems with control of mass:
  a survey, Milan J. Math. 78~(2) (2010) 417--455.


\bibitem{Smith95}
H.~L. Smith, P.~Waltman, The theory of the chemostat, Vol.~13 of Cambridge
  Studies in Mathematical Biology, Cambridge University Press, Cambridge, 1995,
  dynamics of microbial competition.


\bibitem{turing1952chemical}
A.~M. Turing, The chemical basis of morphogenesis, Philosophical Transactions
  of the Royal Society of London B: Biological Sciences 237~(641) (1952)
  37--72.


  

\end{thebibliography}
\end{document}